# Multiobjective Home Appliances Scheduling Considering Customer Thermal Discomfort: A Multistep Look-ahead ADP-Based Approach

Babak Jeddi, Yateendra Mishra, Gerard Ledwich
School of Electrical Engineering and Computer Science
Queensland University of Technology (QUT), Brisbane, Australia

*Abstract*—This paper proposes a multiobjective home energy management unit (MO-HEMU) to balance the electricity payment and thermal discomfort of a household by properly scheduling devices in a time-varying price environment. The thermal discomfort is measured by the deviation of indoor and hot water temperature from the user's ideal temperature. The home devices include solar Photovoltaics, a battery storage system, deferrable, and thermal appliances. The proposed MO-HEMU is formulated as a dynamic program and a method based on the approximated dynamic programming is used as the scheduling algorithm. The developed method, called multistep look-ahead algorithm, is an iterative algorithm that overcomes the *curse of dimensionality* of the exact DP by choosing a decision based on a limited number of stages ahead. The effectiveness of the proposed model is investigated through numerical simulations. The proposed MO-HEMU enables customers to find the desired trade-off between electricity payment and a discomfort level.

*Index Terms*-- Appliances scheduling, approximate dynamic programming, home energy management, multiobjective optimization, thermal discomfort.

## I. Introduction

Demand response (DR), defined as incentivizing customers to alter their electricity usage patterns in response to time-changing prices such as real-time pricing (RTP) [1], is an important aspect of smart networks. Implementation of such a mechanism is facilitated through two-way communication between the customers and utility. DR programs not only benefit the utility and network operator by effectively reducing the network's peak load, but it also enables price-sensitive customers to manage their electricity usage by making appropriate consumption decisions based on the received price signals. In this regard, smart home energy management units (HEMUs) are becoming more popular for residential customers [2]. In simple words, a HEMU optimally schedules the power flows among the devices in a home in response to the price to keep the comfort and economic desire of the inhabitants within a satisfactory level [3, 4]. A HEMU can communicate with the utility to receive the required information, such as weather forecasts and electricity price, and accordingly, find the optimal operation schedule for the devices. The core part of a HEMU is the scheduling algorithm.

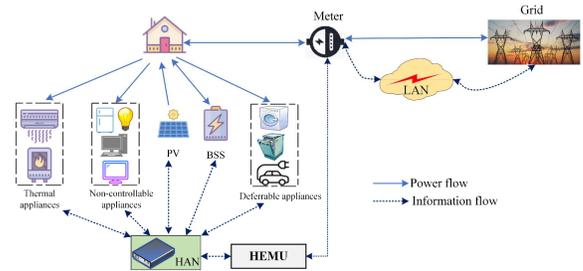

Figure 1. Structure of a smart home with the proposed HEMU

**Related works:** Recently, lots of work conducted on developing efficient scheduling algorithms for HEMUs. A MILP-based approach [5], a model predictive control [6], and PSO algorithm [7] are some examples. A comprehensive review can be found in [2]. When it comes to designing a better HEMU, reaching a balance between optimality and complexity is important. In this regard, dynamic programming (DP) has recently gained more attention [3, 4, 8]. DP is a powerful optimization method, suitable for practical implementations and flexible in terms of handling constraints and objective of any form. However, it suffers the "*curse of dimensionality*" issue, i.e. the computational time increases exponentially as the problem size grows. One solution to this issue is approximate dynamic programming (ADP) [9].

Despite the extensive research, most of the studies consider only a single objective, e.g., minimizing the cost of electricity consumption (CoEC) by scheduling all the appliances to operate at low-price hours which may cause discomfort for the customer. Discomfort refers to the inconvenience experienced by the occupants during load scheduling. A realistic HEMU not only reduces CoEC, but it also minimizes the customer's discomfort at the same time. Simply put, two conflicting and incommensurable objectives should be optimized together leading to a multiobjective HEMU (MO-HEMU). Such a MO-

HEMU, instead of a single solution, provides the customer with a list of electricity payments and the related discomfort levels which are called Pareto solutions. Then, the customer can choose the best-compromise scheduling of appliances from the list based on her own preferences. Some MO scheduling algorithms in [11] and [12] have been proposed in which the MO problem is converted to a single objective one by multiplying each objective to a scalar weight and summing weighted objectives together. However, the proper selection of weights to have a smooth Pareto front is challenging.

**Contributions:** In this study, a MO mixed integer nonlinear programming (MO-MINLP) model for HEMU is formulated considering the thermal discomfort level (TDL) and CoEC. The TDL is measured in terms of the deviation of indoor and hot water temperature from the user's ideal temperatures. The proposed model aims to minimize CoEC and TDL jointly by appropriately determining the schedule of the devices within a home in a time-varying electricity price environment. The home comprises solar PVs, a battery storage system (BSS), deferrable appliances, air conditioning (AC) and electric water heater (EWH). An ADP-based approach called multistep look-ahead algorithm (MLA) is developed as the scheduling algorithm. The MLA is an iterative algorithm that chooses a decision that minimizes the objective function for the current stage plus a limited number of stages ahead, rather than considering all the remaining stages to the end of decision horizon. The MO formulation is based on the $\varepsilon$-constrained technique in which the most preferred objective function (e.g., CoEC) is optimized and the other objective function (e.g., TDL) is considered as an inequality constraint.

**Synopsis:** In Section II, the proposed HEMU model is presented in details. The MO formulation and the MLA are introduced in Section III. Numerical Simulations are given in Section IV followed by conclusions in Section V.

## II. PROPOSED HEMU FORMULATION

The structure of a home equipped with the proposed HEMU is shown in Fig. 1. It receives day-ahead price signals from the utility over a communication infrastructure such as a local area network (LAN). The energy management is done for the day ahead upon available data (i.e., PV generation and household demand) and the resulting schedule is applied to all home devices through a home area network (HAN). The HEMU sells surplus power to the grid at a feed-in-tariff (FiT) price. The decision period is a day including $T=48$ time intervals with the granularity of $\Delta t=30$ minutes.

### A. Model of Controllable Devices

*1) BSS:* The state of charge (SOC) defines the state variable, i.e., $s_t^b \triangleq SOC_t$, where $SOC_t \in [SOC_{min}, SOC_{max}]$. The decision variable $u_t^b$ controls charging ($u_t^b \geq 0$) and discharging ($u_t^b < 0$) constrained by charging/discharging rates, i.e., $u_t^b \in [\Delta SOC_{min}, \Delta SOC_{max}]$. Transition function of BSS state is. The battery power is obtained by:

$$P_t^b = \frac{E^b}{\eta_c^b \Delta t}(SOC_t - \eta_l^b SOC_{t-1})^+ + \frac{E^b \eta_d^b}{\Delta t}(SOC_t - \eta_l^b SOC_{t-1})^- \quad (1)$$

where $(.)^+ = max\{0,.\}$ and $(.)^- = min\{0,.\}$. $E^b$ is the capacity and $\eta_l^b$ models the self-discharging process. $\eta_c^b$ and $\eta_d^b$ are the charging and discharging efficiencies, as shown in Fig. 2 [10].

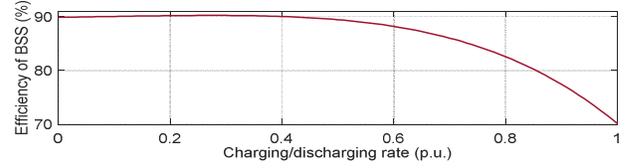
Figure 2. Characteristic of the BSS

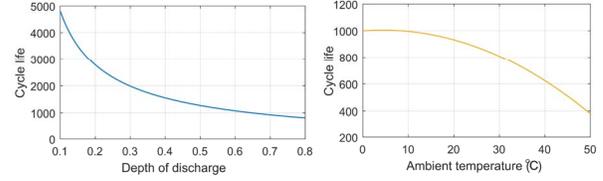
Figure 3. Effect of DoD and temperature on the NiMH battery lifetime

*2) Deferrable non-interruptible appliances:* The starting time of these appliances can be delayed across the day, but once the appliance $a \in A$ is turned on, it has to work for $I_a$ time slots uninterruptedly to finish its task. Dishwasher and laundry appliances belong to this class. The power consumption at each time might be different. Assuming $\boldsymbol{P}_a = \{P_1^a, P_2^a, \ldots, P_{I_a}^a\}$ and $[\alpha_a, \beta_a]$ as the consumption profile and customer-chosen period to activate the appliance $a$, respectively, the goal is to find the best activation time for the appliance. The state is defined by the number of time slots that the appliance has been operated, i.e., $s_t^a \triangleq l_t^a$, where $l_t^a \in \{0,1,\ldots,I_a\}$ is the operational status. Therefore, $s_t^a = 0$ means that the appliance is not activated yet, while $s_t^a = I_a$ indicates that the operation has been completed. The decision variable $u_t^a \in \{0,1\}$ controls the activation of the appliance. The state of appliance depends on the current state and the operating decision. At $t \notin [\alpha_a, \beta_a]$, the appliance must be off (i.e., $s_t^a = 0$) and consequently $u_t^a = 0$. When $t \in [\alpha_a, \beta_a - I_a + 1)$, the controller may turn the appliance on or delay the operation. If the controller does not switch the appliance on, the state remains 0, otherwise the operational status $l_t^a$ increments by 1. The maximum delay time is $t = \beta_a - I_a + 1$ in which the appliance should be activated immediately. When the appliance is activated, it has to work consecutively (i.e., $u_t^a = u_{t+1}^a = \cdots = u_{t+I_a}^a = 1$) until it finishes its task in which case it moves to $s_t^a = I_a$ and has to be idle for the rest of the day (i.e., $u_t^a = u_{t+1}^a = \cdots = u_T^a = 0$).

*3) Deferrable interruptible appliances:* These appliances, such as Plug-in electrical vehicles (PEVs), can be turned on and off multiple times. The controller finds the best time slots within $[\alpha_f, \beta_f]$ for the operation of appliance $f \in F$ with $I_f$ total job length. At $t \notin [\alpha_f, \beta_f]$, $s_t^f = u_t^f = 0$. If the appliance is switched on during $t \in [\alpha_f, \beta_f - I_f + 1)$, the state moves to $s_{t+1}^f = s_t^f + 1$, otherwise $s_{t+1}^f = s_t^f$. At time $t = \beta_f - I_f + 1$, if $s_t^f = 0$ then the appliance should be activated instantly and it has to work consecutively for the next $I_f$ time slots.

*4) Thermal appliances:* The controller adjusts the power consumption of AC and EWH to keep the indoor temperature and temperature of the water inside the tank within the customer desired range. The state variables are defined by the indoor and hot water temperature. That is $s_t^{AC} \triangleq Tem_t^{in}$ and $s_t^{EWH} \triangleq Tem_t^W$ where $Tem_t^{in} \in [Tem_{min}^{in}, Tem_{max}^{in}]$ and $Tem_t^W \in [Tem_{min}^W, Tem_{max}^W]$. The decision variables are $u_t^{AC} \triangleq P_t^{AC}$ and $u_t^{EWH} \triangleq P_t^{EWH}$, where $P_t^{AC} \in \{0, p_1^{AC}, \ldots, p_{I_{AC}}^{AC}\}$ and $P_t^{EWH} \in$

$\{0, p_1^{EWH}, \ldots, p_{I_{EWH}}^{EWH}\}$, which are used to find the best consumption level for AC ($P_t^{AC}$) and EWH ($P_t^{EWH}$) at each time. $I_{AC}$ and $I_{EWH}$ are the total number of consumption levels for AC and EWH. The transition functions are [5, 13]:

$$Tem_t^{in} = \sigma Tem_{t-1}^{in} + (1-\sigma)(Tem_t^{out} - \frac{COP \times P_t^{AC}}{C}) \quad (2)$$

$$Tem_t^W = Tem_{t-1}^W + \frac{(P_t^{EWH} + f_t^w C_p(Tem_{inlet} - Tem_t^W) + A_w(Tem_a - Tem_t^W)) \times \Delta t}{m_w C_p} \quad (3)$$

where $Tem_t^{out}$ (°C), $\sigma$, $C$ and $COP$ are the outdoor temperature, inertia factor, thermal conductivity, coefficient of performance, respectively. $f_t^w$, $C_p$, $Tem_{inlet}$, $Tem_a$, $A_w$, and $m_w$ are hot water consumption, specific heat of water, inlet water temperature, ambient temperature, total thermal conductance, and mass of tank water, respectively. Note that the EWH is located in the area affected by AC operation, so the ambient temperature $Tem_a$ is equal to the indoor temperature $Tem_t^{in}$.

*B. Objective Functions*

The first objective function is the daily CoEC as follows:

$$\min: CoEC = \sum_{t=0}^{t=T-1} \pi_t L_t^+ \Delta t + FiT_t L_t^- \Delta t + C_t^d P_t^b \Delta t \quad (4)$$

where $\pi_t$ and $FiT_t$ are the electricity price and FiT rate, respectively. $L_t$ is the total electricity consumption given by:

$$L_t = \sum_{a \in A} P_{l_t^a}^a u_t^a + \sum_{f \in F} P_{l_t^f}^f u_t^f + P_t^{AC} + P_t^{EWH} + P_t^{un} + P_t^b - P_t^{PV} \quad (5)$$

where $P_t^{PV}$ and $P_t^{un}$ are the PV generation and un-controllable demand, respectively. The positive and negative values for $L_t$ correspond to the consumption and feed-in power cases, respectively. The first and second terms in (4) are the payment and reward for the energy import and export, respectively. The third term is the BSS operational cost in which $C_t^d$ refers to the degradation cost. Ambient temperature and depth of discharge ($DoD=1-SOC$) are the two most dominant factors affecting the lifespan of BSS [14]. The degradation cost is obtained by $C_t^d = C_{BESS}/L_E$ [3], where $C_{BESS}$ and $L_E$ are the replacement cost and the lifetime of BSS. $L_E$ is measured as the number of cycles that it can deliver before it is no longer able to provide sufficient energy [15]. The life cycle is obtained by [14]:

$$L_E = \eta_{Tem} \cdot \eta_{DoD} \cdot L_R \cdot E^b \cdot DoD \quad (6)$$

where $L_R$ is the nominal life cycle in ideal conditions. $\eta_{Tem} = L(Tem)/L_R$ and $\eta_{DoD} = L(DoD)/L_R$ are the temperature and DoD dependence factors, respectively. $L(Tem)$ and $L(DoD)$ are functions relating the battery cycle life for various temperatures and DoDs, obtained from experiment and curve fitting techniques, Fig. 3 [14]. Finally, the degradation cost is:

$$C_t^d = \frac{C_{BESS} \cdot L_R}{L(Tem) \cdot L(DoD) \cdot E^b \cdot DoD} \quad (7)$$

The second objective function is related to the customer's thermal discomfort. Discomfort refers to the inconvenience experienced by the occupants during the load scheduling. The proposed HEMU guarantees to keep the indoor and hot water temperature within the customer's desired range during the day by appropriately scheduling the operation of AC and EWH. However, the middle point of the range is the ideal temperature, $Tem_{idl}$, which can be used to measure the TDL experienced by the occupants. Thus, the discomfort depends on the difference between the scheduled temperature and the ideal temperature. As the difference increases, the TDL increases too. Therefore, the incurred TDL at each time slot $t$ due to the change of indoor temperature is modeled by:

$$TDL_t^{in} = \exp(\frac{|Tem_t^{in} - Tem_{idl}^{in}|}{\Delta Tem_{ther}^{in}}) \quad (8)$$

where $\Delta Tem_{ther}$ is the thermal tolerance for the indoor temperature specified by the user. The TDL can be modeled by other functions and the choice of this exponential function here is based on that it is convex and symmetric around the ideal temperature. This function has a minimum value of 1 when the temperature equals the ideal temperature. As the temperature moves away from the ideal one, more discomfort is experienced by the user. The maximum value is obtained when the maximum deviations from the ideal temperatures happen, i.e., $Tem_t = Tem_{idl} \pm \Delta Tem_{ther}$, which is 2.718. Therefore, at any time slot $t$, $1 \leq TDL_t^{in} \leq 2.718$. The same model formulates the incurred $TDL_t^W$ due to the change in hot water temperature. Finally, the total TDL at time slot $t$ is the summation of individual TDLs incurred by indoor and hot water temperature changes. Then, the daily TDL is obtained:

$$\begin{aligned}\text{Min: } TDL &= \sum_{t=0}^{t=T-1}\{TDL_t^{in} + TDL_t^W\} \\ &= \sum_{t=0}^{t=T-1}\{\exp(\frac{|Tem_t^{in} - Tem_{idl}^{in}|}{\Delta Tem_{ther}^{in}}) + \exp(\frac{|Tem_t^W - Tem_{idl}^W|}{\Delta Tem_{ther}^W})\}\end{aligned} \quad (9)$$

Considering that the entire day is divided into $T=48$ time slots, the daily TDL varies between 96 and 260.95

### III. PROPOSED OPTIMIZATION APPROACH

*A. Exact DP*

The optimal solution for the scheduling problem is found by recursively solving the Bellman equation for t = T-1, ...,0:

$$J_t(\mathbf{s}_t) = \min_{\substack{u_t^a \in \Phi_t^a, u_t^f \in \Psi_t^f, u_t^b \in \Gamma_t^b \\ u_t^{AC} \in M_t^{AC}, u_t^{EWH} \in K_t^{EWH}}} \{g_t(\mathbf{s}_t, \mathbf{u}_t) + J_{t+1}(\mathbf{s}_{t+1}|(\mathbf{s}_t, \mathbf{u}_t))\} \quad (10)$$

$\mathbf{s}_{t+1} = f(\mathbf{s}_t, \mathbf{u}_t)$
$\Gamma_t^b = \{u_t^b | u_t^b \in [\Delta SOC_{min}, \Delta SOC_{max}], SOC_{min} \leq SOC_t \leq SOC_{max}\}$
$K_t^{EWH} = \{u_t^{EWH} | u_t^{EWH} \in \{0, p_1^{EWH}, \ldots, p_{I_{EWH}}^{EWH}\}, Tem_{min}^W \leq Tem_t^W \leq Tem_{max}^W\}$
$M_t^{AC} = \{u_t^{AC} | u_t^{AC} \in \{0, p_1^{AC}, \ldots, p_{I_{AC}}^{AC}\}, Tem_{min}^{in} \leq Tem_t^{in} \leq Tem_{max}^{in}\}$

where $\mathbf{s}_t = (s_t^a, s_t^f, s_t^{AC}, s_t^{EWH}, s_t^b)$ is the total state space and $\mathbf{u}_t = (u_t^a, u_t^f, u_t^{AC}, u_t^{EWH}, u_t^b)$ is the total decision space at time slot $t$. $g_t(s_t, L_t)$ is the one-stage objective function (which can be either of CoEC or TDL) and $J_t(s_t)$ is the optimal cost-to-go at state $\mathbf{s}_t$ and time $t$. $f(.)$ is the transition function of the system that governs the evolution of the state variables over time. $\Phi_t^a, \Psi_t^f, \Gamma_t^b, K_t^{EWH}$ and $M_t^{AC}$ are sets of all feasible operating decisions for each controllable device at stage $t$ based on the related constraints. At the terminal stage $T$, the cost-to-go $J_T(s_T) = 0$. At stage 0, the term $J_0(s_0)$ is the minimized objective function for the entire day for the given initial state $\mathbf{s}_0$. In brief, the exact DP process begins from the end of the decision horizon in which (10) is solved for t = T-1. Then, it steps backward in time and adds one stage at a time, solves a one-stage problem and continues this process until all stages are covered. This approach checks all possible states at any stage $t$ and evaluates all feasible decisions to find the optimal decision which might be computationally intractable.

## B. Multistep Look-ahead ADP

One approach to address the *curse of dimensionality* issue is to settle for a suboptimal solution as a trade-off for a reasonable computational time. The main idea of ADP is to replace the true cost-to-go functions $J_t$ with some approximated functions $\tilde{J}_t$. The most common technique is to start from $t = 1$, step one stage forward in time until the terminal stage is reached, approximate cost-to-go functions for each state being visited, and update approximations iteratively. This approach is called **one-step look-ahead algorithm** (OLA) in which for the state $s_t$ at every stage $t$ only a single stage ahead is observed and the decision that minimizes the following expression is chosen:

$$\hat{u}_t^k = \underset{u_t \in U_t(s_t)}{\arg\min} \{g_t(s_t, u_t) + \gamma \tilde{J}_{t+1}^{k-1}(s_{t+1}|(s_t, u_t))\} \quad (11)$$

where $\tilde{J}_{t+1}^{k-1}$ is an approximation of the true cost-to-go function from previous iterations, $\gamma \in (0,1)$ is the discount factor and $k$ is the iteration number. $U_t(s_t)$ is the decision space based on the feasible decisions for devices operation in (10). Then, the value of being at state $s_t$ and making the decision $\hat{u}_t^k$ is:

$$\hat{v}_t^k = g_t(s_t, \hat{u}_t^k) + \gamma \tilde{J}_{t+1}^{k-1}(s_{t+1}|(s_t, \hat{u}_t^k)) \quad (12)$$

Now, the approximation for the visited state $s_t$ is updated by:

$$\tilde{J}_{t+1}^k(s_t) = (1 - \alpha^k)\tilde{J}_{t+1}^{k-1}(s_t) + \alpha^k \hat{v}_t^k \quad (13)$$

where $\alpha \in (0,1)$ is the stepsize and tunes the convergence speed of the algorithm. The OLA proceeds toward the end of the optimization horizon. Then, it goes for the next iteration and the process starts over again from $t=1$. At each iteration, the cost-to-go approximations of the states being visited are updated. The advantage is that it does not require to loop over all possible states and thus is faster than the DP.

Although OLA looks promising, a large number of iterations may be required to reach a solution close enough to optimal. One way to have a faster convergence and better solution is to use a **two-step look-ahead algorithm** (TLA) in which at stage $t$, the decision is chosen by observing two stages ahead. That is, for every feasible decision $u_t$ that results in $s_{t+1}$, there is:

$$\tilde{J}_{t+1}^k(s_{t+1}) = \min_{u_{t+1} \in U_{t+1}} \{g_{t+1}(s_{t+1}, u_{t+1}) + \gamma \tilde{J}_{t+2}^{k-1}(s_{t+2}|(s_{t+1}, u_{t+1}))\} \quad (14)$$

where $\tilde{J}_{t+2}^{k-1}$ is the approximation of $J_{t+2}$ from previous iterations. The TLA can be extended to a **multi-step look-ahead algorithm** (MLA) with look-ahead of $l > 2$ stages in which case $l$ stages ahead are observed to make the decision.

Naturally, by increasing the size of the look-ahead $l$, the performance is improved, at the expense of increased calculation time. For the OLA, only a single minimization problem in (11) needs to be solved for the visited state $s_t$ per stage. In the TLA, however, this equation has to be solved for each possible state $s_{t+1}$ generated from $s_t$. The computational burden increases for the MLA in the same way. Therefore, a substantial computation might be required when the state or action space is large. One possibility to save computation time is minimization over a subset $\overline{U}_t(s_t) \subset U_t(s_t)$ instead of $U_t(s_t)$.

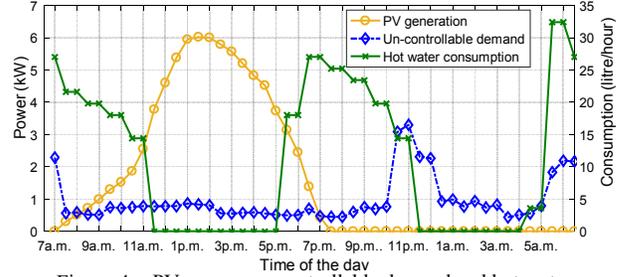

Figure 4. PV power, un-controllable demand and hot water consumption

TABLE I. DEFERRABLE APPLIANCES PARAMETERS

|  | Washing Machine | Dishwasher | PEV |
|---|---|---|---|
| Power (kW) | {0.5, 2, 1.2} | {2.2, 0.3, 1.4, 0.3} | 3.2 |
| Time Interval | 9:00-18:00 | 9:00-18:00 | 18:00-6:00 |
| Required Slots | 3 | 4 | 8 |

The determination of $\overline{U}_t$ is often problem-dependent. One potential option is to select some important variables of the problem and do minimization over the feasible decisions related to these variables, while keeping the decisions of the remaining variables at some nominal values. This approach can be specifically useful in MLA in which the original decision set $U_t(s_t)$ is used for the one-stage look-ahead approximations, while the reduced decision set $\overline{U}_t(s_t)$ is used for the multi-step look-ahead approximations. In the formulated HEMU, for example, the sets of charging/discharging decisions for the BSS, i.e., $r_t^b$, can be nominated to build such a $\overline{U}_t(s_t)$.

## C. Multiobjective Formulation

The proposed MO-HEMU is stated as follows:

$$\begin{array}{ll} \min: & (CoEC(s_t, u_t), TDL(s_t, u_t)) \\ \text{subject to:} & u_t \in U_t(s_t) \end{array} \quad (15)$$

Due to the conflicting nature of the two objectives, a single solution that simultaneously minimizes both objectives is not obtainable. Instead, a set of Pareto solutions can be found. None of the solutions dominates others and each of them is a compromise between the two competing objectives [16].

A well-organized solution method for (15) is the $\varepsilon$-constrained technique in which the CoEC is considered as the main objective function and the TDL is formulated as an inequality constraint [17]. The range of variation for the TDL needs to be determined first. Then, this range is divided into a number of equidistant points Thus, the MO-HEMU in (15) changes to the following form:

$$\begin{array}{ll} \min: & CoEC(s_t, u_t) \\ \text{subject to:} & \begin{cases} TDL(s_t, u_t) \leq \varepsilon_i \\ u_t \in U_t(s_t) \end{cases} \end{array} \quad (16)$$

where

$$TDL^{max} - \left(\frac{TDL^{max} - TDL^{min}}{k_2}\right) \times i, \quad i = 0,1,2,\ldots,k_2 \quad (17)$$

and $TDL^{max}$ and $TDL^{min}$ are the maximum and minimum values of the TDL variation range, which were determined as 96 and 260.95, respectively, in section II.B. In fact, the MO problem is converted into $(k_2 + 1)$ single objective problems each of which with different values of $\varepsilon$ on the right-hand side of the inequality constraint. By solving each of these problems, $(k_2 +$

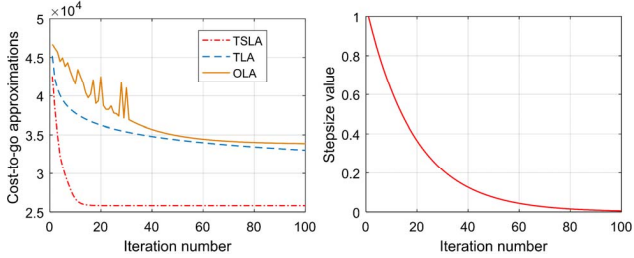

Figure 5. Convergence and stepsize variation rule of ADP algorithms

TABLE II. OBTAINED CoEC (CENTS) BY ALGORITHMS FOR DIFFERENT CASE STUDIES

|  | DP | OLA | TLA | TSLA | | PSO | |
|---|---|---|---|---|---|---|---|
|  | CoEC | CoEC | CoEC | CoEC | Opt.* | CoEC | Opt.* |
| Case 1 | 371.06 | 412.6 | 386.7 | 379.2 | 98% | 371.0 | 100% |
| Case 2 | 442.58 | 521.6 | 464.0 | 484.7 | 91% | 486.5 | 90% |
| Case 3 | 658.91 | 828.6 | 776.1 | 756.5 | 87% | 778.8 | 85% |
| Case 4 | 400.16 | 546.5 | 493.6 | 480.2 | 83% | 491.8 | 81% |
| Case 5 | NA** | 929.7 | 839.7 | 801.2 | NA | 830.3 | NA |

* Optimality, **NA=not available

TABLE III. COMPARISON OF CALCULATION TIME OF ALGORITHMS FOR DIFFERENT CASE STUDIES

|  | $|s_t|$ | $|u_t|$ | Calculation Time (Sec) | | | |
|---|---|---|---|---|---|---|
|  |  |  | DP | OLA | TLA | TSLA |
| Case 1 | 180 | 8 | 2.05 | 1.16 | 1.88 | 2.18 |
| Case 2 | 1620 | 40 | 10.12 | 2.83 | 4.78 | 7.85 |
| Case 3 | 17820 | 80 | 792.64 | 5.67 | 19.94 | 51.50 |
| Case 4 | 11340 | 280 | 1458.3 | 12.64 | 41.17 | 81.64 |
| Case 5 | 124740 | 560 | NA | 189.20 | 288.95 | 375.74 |

* $|s_t|$= number of state variables at each stage, $|u_t|$= number of decision variables at each stage

1) solutions as the Pareto list are found. An advantage of this technique is that it allows the decision maker to have control over the Pareto front by correctly choosing the values of $k_2$. The choice of a higher value for $k_2$ results in more solutions in the Pareto list with better density but at the expense of higher calculation time. In this paper $k_2 = 6$.

## IV. NUMERICAL SIMULATIONS

The proposed MO-HEMU is applied to a typical home with solar PVs and a BSS of size 8 kWp and 5 kWh, respectively. The PV generation, un-controllable demand profile and daily hot water consumption are shown in Fig. 4. The hourly electricity prices and outdoor temperature are given in [18]. The household has three deferrable appliances as listed in Table I. A PEV with a battery capacity of 16 kWh is considered assuming that it arrives at home at 6:00 p.m. with 20% of initial charge and leaves at 6:00 a.m. next day with full energy. Therefore, it requires 12.8 kWh (80%) energy to get fully charged. For the BSS $C_{BESS} = 1500$ \$, $SOC_{min} = 0.2$, $SOC_{max} = 0.8$, $\delta_{SOC} = 0.1$ and $L_R = 3000@85\%\ DoD$ are chosen. The consumption levels of AC and EWH are {0, 1, 1.5, 2, 2.5} and {0, 4} kW, respectively. The user's ideal indoor and hot water temperature are 22℃ and 60℃ with the thermal tolerance of 2 and 5℃, respectively. The step size in discretizing the indoor and hot water temperature are 0.5 and 1℃. AC and EWH parameters are as in [18] and [13]. Also, $FiT$=6 cents/ kWh.

*1) Single objective optimization:* First, in order to evaluate the performance of the developed MLA in terms of the solution quality and computational time, different case studies with a different combination of devices and thus a different number of state and decision variables are defined as follows:

- *Case #1:* only deferrable appliances are considered;
- *Case #2:* deferrable appliances and AC considered;
- *Case #3:* deferrable and thermal appliances considered;
- *Case #4:* deferrable appliances, AC and BSS considered;
- *Case #5:* all devices are included.

A single objective optimization problem for minimizing CoEC is solved for each case study with three MLAs including OLA, TLA, and three-step look-ahead algorithm (TSLA). Moreover, the exact DP is applied to get the optimal solutions as the benchmark. Also, PSO algorithm is applied to each case study [7]. The results and computational time are compared in Table II and III. Please note that it is not computationally feasible to solve case #5 with the exact DP using the available computer. The computational time of the exact DP increases exponentially as the problem dimension enlarges, while the MLAs are much faster than the exact DP and also the increase in calculation time is almost linear. For the MLAs, the performance is improved as the length of look-ahead increases with a slight growth in the computational time. It should also be mentioned that the MLAs perform better when the system dimension is small as they give more close-to-optimal solutions. However, the optimality decreases as the dimensions of the problem grow. Typical convergence curves of three MLAs in solving the same problem are compared in Fig. 5. By increasing the size of the look-ahead $l$ fewer iterations are required to achieve convergence.

One of the challenges for ADP methods is finding the right stepsize $\alpha^k$ in (13) that tunes the convergence of the algorithm [9]. In this paper, a self-adaptive formula is proposed:

$$\alpha^k = (\frac{1}{2 \times iterMax})^{1/iterMax} \times \alpha^{k-1} \quad (18)$$

where iterMax is the maximum number of iterations. This stepsize rule puts the highest weights on the observations of the value of visiting the state $s_t$ on the early iterations. As the algorithm proceeds, the stepsize value drops to zero, and therefore more weights apply on the current approximations $\tilde{J}_{t+1}^{k-1}(s_t)$, as shown in Fig. 5.

*2) Multiobjective optimization:* Finally, the MO-HEMU model is solved by the TSLA considering all devices (case #5) and the obtained Pareto list is displayed in Fig. 6. The choice of the best compromise solution from this list depends on the preferences of the customer; whether minimizing the electricity cost is her priority or minimizing thermal discomfort is more important. To further investigate the differences among the Pareto solutions, three solutions corresponding to the minimum CoEC, average CoEC and TDL, and minimum TDL are chosen. The scheduling of devices and indoor and hot water temperature for each of these solutions are compared in Figs. 7 and 8. As observed, for the case related to the minimum CoEC, the EWH and AC are operated less frequently. Therefore, the largest deviation for the indoor and hot water temperature from the ideal temperature occurs, and thereby the maximum TDL achieved. On the other hand, the EWH and AC are forced to work for a longer time in the case of minimum TDL to keep the indoor and hot water temperature as much as close to the ideal ones, and therefore the minimum deviation happens. The hourly electricity price is also depicted in Fig. 7, which clearly proves that the controller tends to activate home devices during low-price or sunshine hours whenever possible to save the electricity payment. Also, the controller takes the most out of BSS by charging it from the grid at the lowest price and discharging when the morning price peak happens.

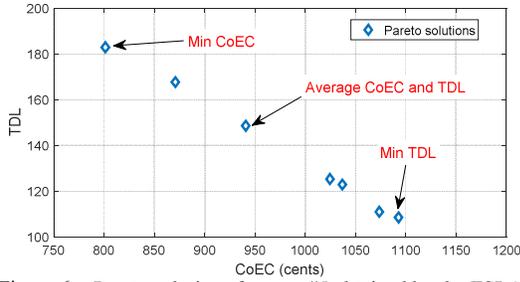

Figure 6. Pareto solutions for case #5 obtained by the TSLA

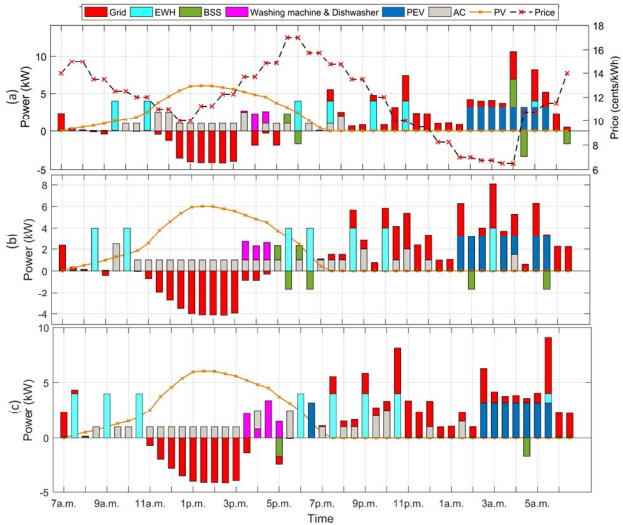

Figure 7. Scheduling of devices corresponding to (a) minimum CoEC, (b) average CoEC and TDL, (c) minimum TDL solution

## V. CONCLUSIONS

This paper formulated a MO-MINLP structure for a HEMU that takes into account the customer's thermal discomfort and electricity payment. The controllable devices of home including thermal and deferrable appliances were modeled in a DP-based setting. Also, solar PVs and BSS were considered. A new ADP-based method called MLA was introduced as the scheduling algorithm. The proposed model aimed to find the appropriate schedule of the devices in order to minimize the two conflicting objective functions, i.e., CoEC and TDL, concurrently. The MO formulation is based on the $\varepsilon$-constrained technique in which the CoEC is optimized as the main objective function and the TDL is considered as a constraint. The formulated HEMU is computationally efficient as the MLA is much faster than the exact DP. Furthermore, the proposed MO-HEMU provides significant benefits to price-sensitive customers enabling them to reach a balance between the electricity payment and a comfortable lifestyle.

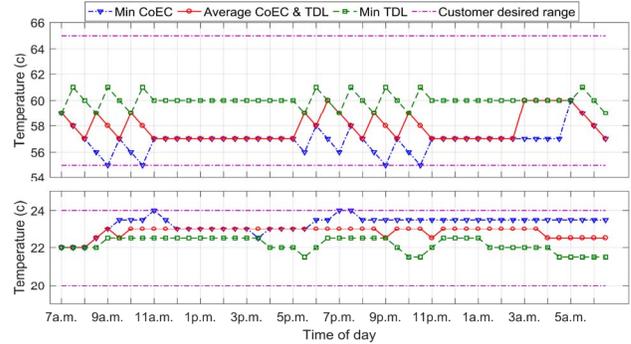

Figure 8. Scheduled hot water and indoor temperature corresponding to the three solutions from the Pareto list